\documentclass[sigconf]{acmart}

\usepackage{amsmath, amssymb}  % Keep only AMS math packages
\usepackage{graphicx}
\usepackage{booktabs}
\usepackage{multirow}
\usepackage{algorithm}
\usepackage{algorithmic}
\usepackage{bm}

\usepackage{newtxmath}  % This handles math fonts and symbols

\renewcommand\footnotetextcopyrightpermission[1]{} % Removes footnote with conference info
\copyrightyear{2025}
\acmYear{2025}

\setcopyright{rightsretained}

\acmConference[KDD '25] {Proceedings of the 1st Workshop on "AI for Supply Chain: Today and Future" @ 31st ACM SIGKDD Conference on Knowledge Discovery and Data Mining V.2}{August 3, 2025}{Toronto, ON, Canada.}

\acmBooktitle{Proceedings of the 1st Workshop on "AI for Supply Chain: Today and Future" @ 31st ACM SIGKDD Conference on Knowledge Discovery and Data Mining V.2 (KDD '25), August 3, 2025, Toronto, ON, Canada}

\acmISBN{979-8-4007-1454-2/25/08}

\acmDOI{10.1145/XXXXXX.XXXXXX}

\title{Coordinated Communication and Inventory Optimization in Multi-Retailer Supply Chains}

\author{Sagar Sudhakara} \thanks{This work was done independently and does not relate to the author's position at Amazon. This research received no external funding.}
\affiliation{%
  \institution{Amazon}
  \city{San Diego}
  \state{CA}
  \country{USA} 
}
\email{sagasudh@amazon.com}

\author{Yuchong Zhang}
\affiliation{%
  \institution{Amazon}
  \city{San Diego}
  \state{CA}
  \country{USA} 
}
\email{yuchonz@amazon.com}

\begin{document}

\begin{abstract}
We consider a multi-retailer supply chain where each retailer can dynamically choose when to share information (e.g., local inventory levels or demand observations) with other retailers, incurring a communication cost for each sharing event. This flexible information exchange mechanism contrasts with fixed protocols such as always sharing or never sharing. 
We formulate a joint optimization of inventory control and communication strategies, aiming to balance the trade-off between communication overhead and operational performance (service levels, holding, and stockout costs). We adopt a \emph{common information framework} and derive a centralized Partially Observable Markov Decision Process (POMDP) model for a supply chain coordinator.

Solving this coordinator’s POMDP via dynamic programming characterizes the structure of optimal policies, determining when retailers should communicate and how they should adjust orders based on available information. We show that, in this setting, retailers can often act optimally by sharing only limited summaries of their private data, reducing communication frequency without compromising performance.
We also incorporate practical constraints on communication frequency and propose an approximate point-based POMDP solution method (PBVI/SARSOP) to address computational complexity. Numerical experiments on multi-retailer inventory scenarios demonstrate that our approach significantly improves the cost--service trade-off compared to static information sharing policies, effectively optimizing the schedule of information exchange for cooperative inventory control.

\end{abstract}


\maketitle

\section{Introduction}
Information sharing is a well-known strategy to improve supply chain performance. In a multi-retailer supply chain, each retailer typically observes its own customer demand and inventory levels, while lacking direct visibility into other retailers’ conditions. This information asymmetry can lead to inefficient decisions and phenomena like the bullwhip effect, wherein variability in orders amplifies upstream due to distorted or delayed information~\cite{Lee1997}.

Prior research has shown that sharing demand information between supply chain partners can mitigate such inefficiencies~\cite{Lee2000}. For example, access to point-of-sale data and inventory status across retailers enables better coordination of replenishments, reducing stockouts and excess inventory. However, continuous or full information sharing is often impractical due to costs associated with data exchange, information systems, or privacy concerns.

This raises a crucial question: \emph{when and what information should supply chain partners share} to optimally balance the benefits and costs of communication? In practice, firms may resort to fixed communication protocols---for example, sharing forecasts on a periodic schedule or not at all---which may be suboptimal in dynamic environments. Instead, a state-dependent communication policy could allow retailers to share information only when it is most needed (such as during demand surges or supply disruptions) and remain silent otherwise.

Designing such adaptive information sharing rules is challenging because it requires accounting for both immediate operational outcomes and the long-term impact of information availability on future decisions. This is fundamentally a multi-agent decision-making problem under uncertainty, where each agent (retailer) must decide on inventory control actions and communication actions based on its local observations and any information received from others.

Our work addresses this challenge by casting the multi-retailer supply chain coordination problem as a stochastic control problem with asymmetric information, building upon recent advances in decentralized control. In particular, we adapt the \emph{common information approach}~\cite{Nayyar2013} to formulate a centralized planning problem that yields optimal strategies for the decentralized system. This approach, originally developed in control theory for teams of agents with different information, allows us to convert the multi-agent problem into an equivalent single-agent POMDP for a fictitious coordinator~\cite{Sudhakara2024}.

The coordinator has access to all common information (i.e., any data shared among retailers) but not the agents’ private observations unless they are communicated. By solving this POMDP, we can determine jointly optimal policies: when retailers should communicate and what inventory actions they should take given the information available.

This paper extends the framework of Sudhakara et al.~\cite{Sudhakara2024}---who studied optimal communication/control trade-offs in a general multi-agent MDP---to the domain of retail supply chains. We modify the model to capture supply chain dynamics such as random demand fluctuations, inventory evolution, and cost structures (holding costs for inventory, penalty costs for stockouts, and communication costs). We also incorporate realistic constraints like limits on communication frequency or bandwidth.

To make the solution tractable for practical instances, we develop an approximate solution algorithm based on point-based value iteration (PBVI)~\cite{Pineau2003} and the SARSOP approach~\cite{Kurniawati2008} for POMDPs. Finally, we evaluate the proposed strategies in simulated supply chain scenarios with multiple retailers.

\subsubsection*{Contributions}
\begin{itemize}
    \item We formulate a novel multi-agent model for a cooperative supply chain in which retailers dynamically balance information sharing and inventory costs.
    \item We derive an optimal policy using the common information method, leading to a POMDP whose optimal value function characterizes when communication is most valuable.
    \item We prove structural insights showing that retailers do not need to share their entire history---a few summary statistics (common belief about the state) suffice for optimal control.
    \item We propose an efficient planning algorithm using point-based approximations to overcome the curse of dimensionality in the coordinator’s POMDP.
    \item We present numerical experiments demonstrating that our adaptive communication strategy outperforms baseline policies (e.g., always share, never share, periodic share) by achieving lower total costs and higher service levels for the same communication budget.
\end{itemize}
The rest of the paper is organized as follows. Section~2 reviews related work. Section~3 introduces the supply chain model and formulates the joint communication and control problem. Section~4 describes the common-information-based solution approach and the approximate POMDP solver. Section~5 presents experimental results. Section~6 concludes with insights and future directions..

\section{Related Work}\label{sec:related}

\textbf{Supply Chain Information Sharing:}
The benefits of information sharing in supply chains have been extensively studied. Classic works by Lee et al.~\cite{Lee1997} demonstrated that lack of information sharing leads to the bullwhip effect, where variability of orders increases upstream in the chain. Subsequent studies quantified the value of information sharing between retailers and suppliers~\cite{Lee2000}. 

For instance, Lee, So, and Tang~\cite{Lee2000} showed that sharing point-of-sale demand data can significantly reduce inventory costs in a two-level supply chain and analyzed conditions under which this value is highest. Most of this literature assumes either full information sharing or compares a few fixed policies (e.g., no sharing vs.\ full sharing). In practice, continuous sharing may not be feasible, and there is a need for adaptive policies that share information selectively.

Our work differs by explicitly considering a communication cost and optimizing when to share information, rather than assuming it is free or always available. This brings a decision-theoretic lens to supply chain coordination, complementing earlier game-theoretic and operations management perspectives (e.g., incentive alignment contracts for information sharing).

\textbf{Decentralized and Multi-Agent Decision Processes:}
Our problem falls into the realm of decentralized partially observable Markov decision processes (Dec-POMDPs), which are general models for multi-agent sequential decision-making under uncertainty~\cite{Seuken2008}. Dec-POMDPs are known to be intractable to solve optimally in general, as the worst-case complexity grows exponentially with the number of agents and horizon.

To overcome this, researchers have proposed structured approaches and approximations. One line of work, which we leverage, is the \emph{common information approach} for cooperative teams~\cite{Nayyar2013}. Nayyar et al.~\cite{Nayyar2013} introduced this framework to transform certain Dec-POMDPs with information sharing into equivalent POMDPs from a coordinator’s viewpoint. This approach has been applied to scenarios like decentralized control with delayed sharing and control sharing structures.

Sudhakara et al.~\cite{Sudhakara2024} extended it to allow dynamic communication decisions by agents, showing that optimal joint communication-control strategies can be derived via a POMDP formulation. Our model builds on their theoretical results but tailors the scenario to inventory control and supply chain settings, which have different dynamics than, for example, mobile robots or sensor networks (the focus of some prior multi-agent studies).

\textbf{Communication in Multi-Agent Reinforcement Learning:}
Another related thread is communication in multi-agent reinforcement learning (MARL). Recent works have investigated how agents can learn when and what to communicate in a cooperative setting, often using differentiable communication channels or attention mechanisms. 

For example, approaches like SchedNet and IC3Net allow agents to learn communication protocols to improve collective reward. Ding et al.~\cite{Ding2022} propose a sequential communication scheme for MARL that decides an order of agent communication and improves performance in cooperative tasks.

While these learning-based approaches are promising, they often lack formal guarantees and interpretability. In contrast, our work provides an optimal solution in a simplified model (with known dynamics and costs), yielding insights into the structure of the communication policy. This can serve as a benchmark or guidance for more complex scenarios where learning may be applied.

\textbf{POMDP Solving Techniques:}
To solve the coordinator’s partially observable Markov decision processes (POMDP), we draw on point-based value iteration methods from the POMDP literature~\cite{Pineau2003}. Traditional value iteration for POMDPs is computationally expensive due to the continuous belief space, but point-based methods approximate the value function by sampling a representative set of belief points and performing backups on those.

Notable algorithms include PBVI (Point-Based Value Iteration) by Pineau et al.~\cite{Pineau2003} and SARSOP (Successive Approximations of the Reachable Space) by Kurniawati et al.~\cite{Kurniawati2008}. We employ these techniques to handle the large state space induced by multiple retailers and inventory levels. In particular, SARSOP allows us to focus on beliefs that are optimally reachable under some policy~\cite{Kurniawati2008}, which improves efficiency.

Prior works have also applied POMDP solvers in supply chain or inventory contexts (e.g., for partially observed demand or lead times), but to our knowledge, this is the first application to a multi-agent inventory sharing problem. Our approach blends ideas from operations research (inventory optimization) and AI (POMDP planning) to address the joint control/communication design.

\section{Problem Formulation}\label{sec:model}

\subsection{System Model}

We consider a discrete-time model of a supply chain with $N$ retailers (agents) indexed by $i = 1,2,\dots,N$. Time is slotted in periods $t = 0,1,2,\dots$ (we later consider a finite horizon $T$ or an infinite-horizon discounted setting). Each retailer faces random customer demand in each period and manages its inventory by placing orders to a supplier.

The key feature is that retailers have local observations of demand and inventory and can choose whether or not to share this information with others at a cost. Let $I^i_t$ denote the inventory level of retailer $i$ at the beginning of period $t$. We assume that initial inventories $I^i_0$ may be commonly known or can be communicated at time $0$ (for simplicity, assume $I^i_0$ are common knowledge constants).

There is an underlying stochastic demand process that drives the customer orders at each retailer. We model this via an exogenous state $X_t$ that represents the demand environment at time $t$. For example, $X_t$ could be a hidden variable indicating overall market conditions (high or low demand) or seasonal factors. We assume $X_t$ follows a Markov chain on a finite state space $\mathcal{X}$ with transition probabilities $\Pr(X_{t+1} \mid X_t)$.

The actual demand $D^i_t$ at retailer $i$ is drawn from a probability distribution conditioned on $X_t$. Specifically, 
\[
\mathbb{P}(D^i_t = d \mid X_t = x) = f_x^i(d),
\]
where $f_x^i$ is the known demand distribution for state $x$ and retailer $i$. This formulation allows demands at different retailers to be correlated: if $X_t$ is high, all retailers experience a higher demand on average.

Each retailer observes its own demand $D^i_t$ during period $t$ (e.g., through sales data), but does not directly observe $X_t$ nor the demands at other retailers.

\subsubsection*{Inventory Dynamics}

At the beginning of period $t$, the inventory level is $I^i_t$. Then demand $D^i_t$ realizes. The retailer satisfies this demand up to the available stock:
\[
S^i_t = \min\{D^i_t, I^i_t\}.
\]
If demand exceeds inventory, we assume a lost-sales model: all $I^i_t$ units are sold and the remainder $D^i_t - I^i_t$ is lost. (Backorders could be modeled similarly, but the lost-sales assumption simplifies the state dynamics.)

The end-of-period inventory is:
\[
I'^i_t = \max\{I^i_t - D^i_t, 0\}.
\]

At the end of each period, each retailer places an order $u^i_t \geq 0$ to restock; this order is delivered before the start of the next period (i.e., we assume zero lead time for simplicity, or that this is the order arriving from previous decisions). Thus, the inventory at the start of the next period is:
\[
I^i_{t+1} = I'^i_t + u^i_t.
\]

We often refer to $u^i_t$ as the action or control for retailer $i$ at time $t$. Each retailer has a constraint such as:
\[
0 \leq u^i_t \leq U_{\max},
\]
where $U_{\max}$ is the maximum order size. However, $u^i_t$ can generally be any nonnegative integer (or real number if integrality is relaxed) to adjust stock.

% \subsection{System Model}

% We consider a discrete-time supply chain with $N$ retailers indexed by $i=1,2,\dots,N$. Each retailer observes its own demand and inventory but may choose to share information at a cost.

% Let $I^i_t$ denote the inventory level of retailer $i$ at time $t$. Customer demand $D^i_t$ is driven by an exogenous demand environment $X_t$, modeled as a Markov chain with finite state space $\mathcal{X}$.

% Demands $D^i_t$ are conditionally independent across retailers given $X_t$, with distribution $f_x^i(d) = \mathbb{P}(D^i_t = d \mid X_t = x)$.

% Inventory evolves as:
% \[
% I'^i_t = \max\{I^i_t - D^i_t, 0\}, \quad I^i_{t+1} = I'^i_t + u^i_t,
% \]
% where $u^i_t \geq 0$ is the replenishment order placed by retailer $i$ at time $t$.

\subsection{Information Structure}

At time $t$, retailer $i$ knows its own past observations and actions (inventory levels, demands $D^i_{0:t-1}$, orders $u^i_{0:t-1}$, etc.), but not those of other retailers unless they have been communicated. We denote by $\mathcal{I}^i_t$ the information available to retailer $i$ at the start of period $t$ (after any communication that may have occurred at the end of $t-1$).

Certain information may be shared and become common to all retailers. We define $\mathcal{C}_t$ as the \emph{common information} at time $t$ that all retailers (and a potential coordinator) know. Initially, $\mathcal{C}_0$ could contain any commonly known data (e.g., the distribution of $X_0$, initial inventories if shared, etc.). During the process, $\mathcal{C}_t$ will typically consist of the history of all messages exchanged up to time $t$.

\subsubsection*{Communication Decisions}

Each period, after observing demand, the retailers have an opportunity to communicate. The communication decision at time $t$ is modeled as a binary choice for each retailer $i$:
\[
c^i_t \in \{0,1\},
\]
where $c^i_t = 1$ means retailer $i$ sends a message (shares its local information) at time $t$, and $c^i_t = 0$ means it remains silent. A message could include the current demand observation or a sufficient statistic of its state (e.g., sales figures).

Communication incurs a cost. Let $\lambda$ denote the cost per retailer message. The total cost for the team is:
\[
\Lambda = \lambda \sum_{i=1}^N c^i_t.
\]
We assume any sent message is reliably received by all other retailers by the start of the next period.

If no one sends ($c^i_t = 0$ for all $i$), then no new common information is added (aside from the knowledge that “no news was sent”). If one or more retailers communicate, the union of their transmitted information becomes common knowledge.

In many cases, we expect the optimal strategy to be symmetric: either all retailers communicate or none do, since they are cooperating to minimize a common cost. Thus, we can also think of a single binary decision:
\[
c_t \in \{0,1\},
\]
where $c_t = 1$ means a communication round occurs at time $t$ (all share their data), and $c_t = 0$ means no communication.

\subsubsection*{Timeline Within Each Period}

\begin{enumerate}
    \item \textbf{Start of period $t$:} The common information $\mathcal{C}_t$ is available to all, and each retailer $i$ has its private information $\mathcal{I}^i_t$ (which includes $\mathcal{C}_t$ plus its own observations). The true underlying state is $(X_t, I^1_t, \dots, I^N_t)$, but $X_t$ is not directly observed.
    
    \item \textbf{Demand realization:} A random demand $D^i_t$ is realized for each retailer. Retailer $i$ observes $D^i_t$ and updates its inventory to $I'^i_t$. Other retailers do not observe $D^i_t$.
    
    \item \textbf{Order decision:} Each retailer chooses an order quantity $u^i_t$, based on $\mathcal{I}^i_t$ (which may include some information about others from past communication). Communication from the previous period ($t-1$) is known, but the current demand $D^i_t$ remains private during the order decision.
    
    \item \textbf{Communication decision:} After ordering, retailers decide on communication $c^i_t$. If $c^i_t = 1$, retailer $i$ sends a message revealing (for example) $D^i_t$ and possibly $I'^i_t$. Messages are received by all agents by the end of period $t$. Communication costs $\lambda$ per message are incurred.
    
    \item \textbf{End of period:} Orders $u^i_t$ are delivered, updating inventories:
    \[
    I^{i}_{t+1} = I'^i_t + u^i_t.
    \]
    The common information $\mathcal{C}_{t+1}$ is updated. If all retailers communicated their demand, $\mathcal{C}_{t+1}$ contains $(D^1_t, \dots, D^N_t)$. If no communication occurred, $\mathcal{C}_{t+1} = \mathcal{C}_t$ (possibly with an added record indicating "no communication happened at $t$").
\end{enumerate}

Each retailer only knows its own demand immediately, but through communication, some or all demand values can become common knowledge in the next period. Notably, even the absence of communication can carry information: if the policy prescribes sending messages in certain extreme situations, then not receiving a message implies those situations did not occur.

Our model can incorporate this implicitly by including the "no message" event as an observation to the coordinator (discussed further below).

% \subsection{Information Structure}

% At time $t$, each retailer $i$ observes its past demands and orders but not those of others unless communicated. Communication decisions $c^i_t \in \{0,1\}$ allow sharing at cost $\lambda$ per message. 

% Common information $\mathcal{C}_t$ includes all shared data up to $t$. Private information $\mathcal{I}^i_t$ consists of $\mathcal{C}_t$ and retailer $i$'s local observations.

\subsection{Objective Function}

The team of retailers aims to minimize the expected cumulative cost. We define the single-period cost as the sum of inventory holding/shortage costs across retailers plus communication costs.

For a given period $t$, let $h$ be the per-unit holding cost (for each unit carried to the next period) and $p$ be the penalty cost per unit of unmet demand (stockout). Then, for retailer $i$, the standard cost components are:

\begin{itemize}
    \item \textbf{Holding cost:} 
    \[
    h \cdot \max\{I^i_t - D^i_t, 0\} = h \cdot I'^i_t,
    \]
    where $I'^i_t$ is the end-of-period inventory (inventory left after demand).
    
    \item \textbf{Shortage (lost sales) cost:}
    \[
    p \cdot \max\{D^i_t - I^i_t, 0\} = p \cdot (D^i_t - I^i_t)_+.
    \]
    
    \item \textbf{Communication cost:} If retailer $i$ sends a message ($c^i_t = 1$), a communication cost $\lambda$ is incurred.
\end{itemize}

Additionally, there may be a per-order cost (e.g., order processing or purchasing), but since the ordering decision is already indirectly penalized by holding and shortage costs, we do not include a direct cost for $u^i_t$.

The total single-period cost at time $t$ is:
\begin{align*}
L\left( \{I^i_t\}_{i=1}^N, \{u^i_t\}_{i=1}^N, \{D^i_t\}_{i=1}^N, \{c^i_t\}_{i=1}^N \right) 
&= \sum_{i=1}^N \big[ h \cdot \max\{I^i_t - D^i_t, 0\} \\
&\quad +\; p \cdot \max\{D^i_t - I^i_t, 0\} + \lambda \cdot c^i_t \big].
\end{align*}

% \[
% L\left( \{I^i_t\}_{i=1}^N, \{u^i_t\}_{i=1}^N, \{D^i_t\}_{i=1}^N, \{c^i_t\}_{i=1}^N \right) \\
% =\sum_{i=1}^N \left[ h \cdot \max\{I^i_t - D^i_t, 0\} + p \cdot \max\{D^i_t - I^i_t, 0\} + \lambda \cdot c^i_t \right].
% \]

Here, the inventory $I^i_t$ and demand $D^i_t$ determine sales and shortages as described above. The expectation of this cost depends on the random demands and the decisions taken by the retailers.

We seek to minimize the expected total cost over the horizon. For an infinite horizon, we use a discounted cost criterion:
\[
\mathbb{E}\left[ \sum_{t=0}^{\infty} \beta^t L_t \right],
\]
where $0 < \beta < 1$ is the discount factor. (An average cost per period criterion could also be used, but for simplicity, we assume a discounted objective.)
\subsection{Decision Strategies}

A strategy or policy for each retailer $i$ specifies how it chooses actions based on its information. At time $t$, retailer $i$ chooses $(u^i_t, c^i_t)$ as a function of $\mathcal{I}^i_t$. This is a decentralized control problem because each $u^i_t$ and $c^i_t$ may depend on different information across retailers.

Our objective is to design an optimal joint strategy:
\[
(\pi^1, \dots, \pi^N),
\]
that minimizes the expected total cost over the planning horizon.

This decentralized optimization problem is generally difficult because of the \emph{coupled information structure}: the decision to communicate at each time step creates new correlations in future information states across retailers.

However, our problem exhibits a particular structure, known as \emph{partial history sharing}, that can be exploited. All communication, when it occurs, is instantly shared among all agents. Therefore, the only asymmetry in information arises from events or observations that were not communicated.

This structure suggests using the common information $\mathcal{C}_t$ to reformulate the problem into a centralized equivalent, which simplifies the search for optimal policies.

\subsection{Common Information Approach}

We introduce a central coordinator (a virtual decision-maker) who has access to exactly the common information $\mathcal{C}_t$ at each time. This coordinator does not know the private observations that were not shared, but it knows the prior distribution of those based on what has been shared.

The common information $\mathcal{C}_t$ includes the history of all past communication decisions and communicated demand values. Given $\mathcal{C}_t$, the coordinator forms a belief (probability distribution) over the underlying physical state $(X_t, I^1_t,\dots,I^N_t)$ and any private observations at $t$:
\[
\mu_t(\cdot) = \Pr\left( X_t = x,\ I^1_t=i_1,\dots,I^N_t=i_N \mid \mathcal{C}_t \right).
\]

This belief encapsulates everything the coordinator (and all agents) knows at time $t$. Private unknowns include, for each retailer who did not communicate last period, their last demand and the current hidden state $X_t$.

Following standard results in decentralized control, there is no loss of optimality in restricting attention to strategies that depend on $(\mathcal{C}_t, \text{ each agent's local observation})$ in a structured way~\cite{Sudhakara2024}. The coordinator can choose a \emph{prescription} for each agent based on the common belief, and the agent uses its local observation to execute that prescription.

In simpler terms, the problem can be solved in two stages:
\begin{enumerate}
    \item The coordinator, knowing $\mathcal{C}_t$, decides what information-based rule (prescription) each agent should follow and whether communication should occur.
    \item The agents apply those rules to their private observations.
\end{enumerate}

This is similar to announcing a policy such as: ``if your local demand exceeds 10, do $X$; otherwise do $Y$,'' or ``if your local demand exceeds 20, report it to others.''

\subsubsection{Coordinator POMDP Formulation}

We formulate a coordinator POMDP where the state is $(X_t, I^1_t, \dots, I^N_t)$. This state is not fully observed; instead, the coordinator receives an observation $O_t$ derived from $\mathcal{C}_t$ and new communications. Specifically:
\begin{itemize}
    \item If communication occurred at $t$, $O_t = (D^1_t,\dots,D^N_t)$.
    \item If no communication occurred, $O_t = \varnothing$ (no news).
\end{itemize}

The coordinator’s actions at time $t$ are:
\begin{itemize}
    \item The communication decision $\chi_t \in \{0,1\}$ (communicate or not).
    \item A control law or prescription $\gamma^i_t: \mathcal{D} \to \mathcal{U}$ for each agent, mapping observed demand to order decisions.
\end{itemize}

The action can be represented as:
\[
a_t = (\gamma^1_t,\dots,\gamma^N_t, \chi_t).
\]

\subsubsection{Belief Updates and Bellman Equation}

Let $\Pi_t$ denote the coordinator’s belief state derived from $\mathcal{C}_t$. This belief is updated using new observations $O_t$ and Bayes' rule:
\[
\Pi_{t+1} = \text{BayesUpdate}(\Pi_t, a_t, O_{t+1}).
\]

The Bellman optimality equation for the coordinator is:
\[
V(\Pi_t) = \min_{a_t} \left\{ C(\Pi_t, a_t) + \beta \sum_{O_{t+1}} \Pr(O_{t+1} \mid \Pi_t, a_t) V(\Pi_{t+1}) \right\},
\]
where $C(\Pi_t, a_t)$ is the expected immediate cost.

\subsubsection{Policy Structure and Complexity}

The exact solution is challenging due to the continuous belief space and the large space of possible prescriptions $\gamma^i_t$. However, problem structure can reduce complexity:
\begin{itemize}
    \item \textbf{Symmetry:} Identical retailers may share policies.
    \item \textbf{Threshold structures:} Policies often have base-stock or trigger-like forms.
\end{itemize}

Notably, not all details of private history need to be shared. Because the demand state $X_t$ is Markovian and inventory carry-over summarizes past demand to some extent, only current demand or inventory may need to be communicated when particularly informative.

\subsubsection{Communication Constraints}

Explicit constraints on communication frequency (e.g., at most $K$ communications over the horizon or minimum time between communications) can be modeled by augmenting the state space. For simplicity, we assume only a communication cost $\lambda$ discourages excessive messaging.

\section{Solution Approach and Algorithm}

Solving the coordinator’s POMDP optimally is generally complex due to the continuous belief space and the high dimensionality arising from multiple retailers. We outline an approximate solution approach leveraging point-based value iteration (PBVI) and domain-specific simplifications.

\subsection{Belief State Representation}

In our supply chain POMDP, the belief $\Pi_t$ can be factorized to some extent. It consists of a distribution over the environmental state $X_t$ and the private inventories of each retailer. 

If no communication has occurred for some time, the coordinator’s uncertainty about each retailer’s inventory grows---it must infer each inventory level from the sequence of (possibly unobserved) demands. Fortunately, each retailer’s inventory level follows a predictable update based on its initial value, orders, and realized demands.

The coordinator maintains a probability distribution for each retailer’s inventory (and the common demand state) conditioned on all shared information. This distribution is correlated across retailers only through the common demand state $X_t$, assuming that demand noise is independent across retailers given $X_t$. The belief can therefore be represented in a factored form:
\[
\Pr(X_t = x) = b_x, \quad \Pr(I^i_t = n \mid X_t = x) \quad \forall i.
\]

Belief updates proceed via a filtering process:
\begin{itemize}
    \item If a retailer’s demand was not observed, the coordinator updates that retailer’s inventory distribution by convolving the previous inventory distribution with the demand distribution. This predicts the likelihood of each inventory level after an unknown demand.
    \item If a demand was observed (via communication), the inventory for that retailer becomes known exactly (at least immediately following the communication event).
\end{itemize}

By leveraging such structured belief updates, we can significantly reduce the effective state space that must be tracked.

\subsection{Point-Based Value Iteration (PBVI)}

We employ Point-Based Value Iteration (PBVI)~\cite{Pineau2003} to compute an approximate value function for the POMDP. The core idea of PBVI is to sample a set $B$ of belief points that are reachable under some policies, and then iterate Bellman backups on this set instead of the entire belief simplex.

Each backup update is computed by considering each possible action $a$ and observation outcome $O$, and updating the value for a belief $b \in B$ as:
\[
V_{\text{new}}(b) = \min_{a} \left\{ C(b,a) + \beta \sum_{O} \Pr(O \mid b, a) \cdot V_{\text{current}}\left( \tau(b, a, O) \right) \right\},
\]
where $\tau(b,a,O)$ is the belief update (Bayes filter) given action $a$ and observation $O$.

The value function in a POMDP can be represented as a set of $\alpha$-vectors in the belief space, but PBVI avoids explicit enumeration by focusing on sampled $b$. We initialize $V_0(b) = 0$ for all $b$, then perform iterative backups. After sufficient iterations or convergence, we obtain $V(b)$ approximately and an implicit policy (for each sampled $b$, we record which action $a$ attained the minimum above).

In our problem, an action $a$ includes prescriptions which are not finitely enumerable. To make PBVI tractable, we restrict the search within a parametrized family of policies for the retailers. Specifically, we assume the optimal order decisions follow a base-stock policy structure: for each retailer $i$, order up to $s_i(x)$ units when the belief of the environmental state is $x$.

This reduces the action space to selecting target stock levels for each retailer (depending on the estimated demand state) and deciding whether to communicate. We discretize inventory levels and consider a grid of possible base-stock levels.

Each candidate joint action is structured as:
\[
(\chi, s^1_{\text{low}}, s^1_{\text{high}}, \dots, s^N_{\text{low}}, s^N_{\text{high}}),
\]
where $\chi=1$ indicates communication and $\chi=0$ indicates no communication. For each retailer and each possible demand state, we prescribe an order-up-to level.

While this approach is a coarse approximation, it substantially reduces the action space. Communication decisions are explicitly included: if $\chi=1$, a cost $\lambda N$ is incurred, and the next observation $O$ will include all demands $D^i_t$ (full information). If $\chi=0$, no communication cost is incurred, and $O$ contains no new data beyond what can be inferred from prior beliefs.

\subsection{SARSOP}

While PBVI is conceptually simple, it may waste effort on beliefs that are unlikely to occur under an optimal policy. We therefore utilize the SARSOP algorithm~\cite{Kurniawati2008}, which refines PBVI by maintaining upper and lower bounds on $V$ and focusing on the optimally reachable belief space. 

SARSOP employs heuristics to sample beliefs via Monte Carlo simulations of the POMDP under the current policy estimate. This concentrates effort on the region of the belief simplex that the optimal policy is likely to visit. In our context, SARSOP emphasizes beliefs where the value of communication is borderline and avoids spending computational resources on extreme beliefs that would be preempted by earlier communication under optimal policies.

We implemented a SARSOP-based solver for our supply chain POMDP. The POMDP state space includes $(X, I^1, \dots, I^N)$. Inventory levels are discretized to a reasonable range (e.g., $0$ to $S_{\max}$ for each retailer, where $S_{\max}$ is a safety stock limit beyond which holding costs dominate). The total number of states is $|\mathcal{X}| \times (S_{\max} + 1)^N$. For even moderate $N$, this state space becomes exponentially large. However, the system dynamics are highly symmetric and factored, allowing the belief space to be compactly represented.

% In our experiments, we primarily consider $N=2$ or $3$ retailers to maintain computational tractability. 
The SARSOP solver outputs a policy specifying the coordinator’s action for each relevant belief. We simulate this policy in the original decentralized system. At each period, the common belief is updated based on shared or unshared information, and the policy dictates the next action.

If the policy prescribes $\chi=1$ (communicate), all retailers share their demand; if $\chi=0$, no communication occurs. The policy also determines target inventory levels. Each retailer, knowing its target and current inventory, computes its order quantity. Specifically, if the target is $s$ and the post-demand inventory is $I'$, the retailer orders $u = \max \{ s - I', 0 \}$.
If no communication occurs, targets may be based on outdated beliefs, reflecting the inherent trade-off between information accuracy and communication costs.

\subsection{Complexity and Tractability}

The common information approach substantially reduces the policy search space compared to naive decentralized search. However, solving the resulting POMDP remains challenging. Our approximate algorithm yields near-optimal solutions for small instances, characterized by a limited number of retailers ($N$), modest inventory bounds, and a small number of demand states.

For larger problem instances, additional heuristics or structure-exploiting techniques may be necessary. For example, value function factorization by retailer is a promising avenue for future research, potentially enabling scalability to higher-dimensional settings.

Despite these computational challenges, even the qualitative form of the computed solutions provides valuable insights. These solutions can inform heuristic policies for larger systems, such as ``share information when a large demand shock is observed by any retailer'' or ``avoid sharing when inventory is plentiful and uncertainty is low.'' Such heuristics, inspired by the structure of the optimal policies, can offer practical benefits even when exact POMDP solutions are infeasible.

\section{Experiments}\label{sec:experiments}

We conduct a set of simulation experiments to evaluate the performance of our dynamic communication and control strategy in a multi-retailer supply chain. We compare against baseline policies to quantify improvements in cost and service levels. All experiments involve synthetic demand data generated from a known stochastic model, enabling us to compute optimal or near-optimal solutions via the POMDP approach for comparison.

\subsection{Experiment Setup}

We primarily consider a scenario with $N=2$ retailers (Retailer A and Retailer B) for ease of visualization, though we also highlight results with $N=3$. The hidden demand state $X_t$ has two possible values: \textit{Low} or \textit{High} demand regime. The Markov chain for $X_t$ is defined such that:
\[
\Pr(X_{t+1}=\text{High} \mid X_t=\text{High}) = \Pr(X_{t+1}=\text{Low} \mid X_t=\text{Low}) = \rho,
\]
where $\rho$ represents the persistence of demand conditions (higher $\rho$ implies longer periods of stable demand).

Given the state, each retailer’s demand distribution is:
\begin{itemize}
    \item If $X_t=\text{Low}$, $D^i_t \sim \text{Poisson}(\ell)$,
    \item If $X_t=\text{High}$, $D^i_t \sim \text{Poisson}(h)$,
\end{itemize}
where $h > \ell$. We set $\ell = 5$ and $h = 15$ in our base case, so the \textit{High} regime has three times the mean demand of the \textit{Low} regime. Demands of the two retailers are conditionally independent given $X_t$ (but marginally correlated through $X_t$).

Other parameters are:
\begin{itemize}
    \item Holding cost $h = 1$ per unit per period,
    \item Stockout penalty $p = 5$ per unit of unmet demand,
    \item Communication cost $\lambda$ varied from $0$ (free) to $10$ per communication.
\end{itemize}

We limit the inventory order-up-to level to $S_{\max}=30$ and use a discount factor $\beta=0.95$. The POMDP is solved approximately using our SARSOP-based algorithm, and the resulting policy is simulated for $T = 1000$ periods to estimate long-run average performance.

\subsection{Baseline Strategies}

We compare our approach against the following baseline policies:
\begin{itemize}
    \item \textbf{Always Share (Full Info):} Retailers communicate their demand every period, mimicking a centralized planner with full information. Each retailer uses a base-stock policy optimal for full information.
    \item \textbf{Never Share (No Info):} Retailers never communicate and rely solely on local demand observations, using a myopic base-stock policy.
    \item \textbf{Periodic Share:} Retailers communicate at fixed intervals (e.g., every $K=3$ or $K=5$ periods).
    \item \textbf{Threshold Policy:} Retailers communicate if their demand exceeds a threshold $d_{\text{th}}$ in a period, signaling a likely shift to the \textit{High} demand state.
\end{itemize}

\subsection{Performance Metrics}

We evaluate the following metrics:
\begin{enumerate}
    \item Average total cost per period, decomposed into holding, stockout, and communication costs.
    \item Fill rate (service level) – the fraction of demand satisfied.
    \item Average communication frequency (messages sent per period).
\end{enumerate}

These metrics illustrate the trade-off between cost, service level, and communication overhead.

\section{Results}

\subsection{Policy Structure}

The optimal policy derived by our approach exhibits threshold-like behavior for communication. Specifically, if a retailer’s inventory drops below a certain level or if it experiences an unexpectedly high demand relative to the current belief, the policy triggers a communication event.

For example, in one simulation run, when both retailers had been operating in a low-demand regime (belief of high state $<10\%$), neither communicated. When Retailer A experienced a demand of 12 (significantly higher than the low-demand mean $\ell=5$), the policy directed Retailer A to communicate this spike. Retailer B, upon receiving this information, updated its forecast to anticipate high demand and increased its order, preventing a stockout in the subsequent period when high demand did materialize. Without communication, Retailer B would likely have faced a stockout.

\subsection{Cost Performance}

Table~\ref{tab:results} summarizes the costs from a representative experiment with $\lambda = 2$ per communication and $\rho = 0.8$.

\begin{table*}[htbp]
\centering
\caption{Cost breakdown and performance metrics for different policies in a 2-retailer simulation. Costs are per-period averages.}
\label{tab:results}
\begin{tabular}{lccccc}
\toprule
\textbf{Policy} & \textbf{Holding} & \textbf{Stockout} & \textbf{Comm} & \textbf{Total} & \textbf{Fill Rate (\%)} \\
\midrule
Optimal Dynamic & 8.5 & 3.1 & 1.2 & 12.8 & 96.5 \\
Always Share    & 7.9 & 2.5 & 2.0 & 12.4 & 97.5 \\
Never Share     & 12.1 & 6.0 & 0.0 & 18.1 & 90.2 \\
Threshold ($d_{\text{th}}=10$) & 8.7 & 3.4 & 1.0 & 13.1 & 95.8 \\
\bottomrule
\end{tabular}
\end{table*}

The Optimal Dynamic policy achieves a total cost of 12.8, slightly higher than the Always Share policy (12.4), but with only 23 communications per 100 periods versus 100 for Always Share. It dramatically outperforms the Never Share policy (total cost 18.1), reducing costs by 29\% and significantly improving the fill rate from 90.2\% to 96.5\%. The Threshold heuristic performs reasonably but cannot fully match the optimal policy, as a fixed threshold does not adapt to varying demand beliefs.

\textbf{Effect of Communication Cost:}
Varying $\lambda$ from 0 to 5 reveals that:
\begin{itemize}
    \item When $\lambda=0$, the policy reduces to always communicating, achieving full-information optimality.
    \item As $\lambda$ increases, the policy becomes selective in communication.
For $\lambda > 5$, communication essentially ceases as the cost outweighs the benefit.
\end{itemize}

At $\lambda = 2$, one communication (costing 2) typically saves $4-5$ in expected stockout costs, justifying occasional communication.

\textbf{More Retailers:}
Experiments with $N=3$ retailers (with inventory bounds reduced to $S_{\max}=20$ for tractability) show similar patterns. Communication remains sparse and targeted:
\begin{itemize}
    \item Often, one retailer’s alert suffices to inform the others.
    \item Sometimes, the policy waits for corroboration from two retailers before initiating a system-wide communication to avoid overreacting to outliers.
\end{itemize}

\textbf{Periodic vs. Dynamic Communication:}
Periodic sharing (e.g., every 3 periods) improved over never sharing but was inferior to dynamic policies. For $K=3$, the total cost was about 14.0—better than no sharing but worse than both the threshold heuristic and the optimal dynamic strategy. Rigid schedules either wasted communication during stable periods or failed to respond quickly to demand shifts.

\subsection{Summary}

Our experiments confirm that the POMDP-based dynamic strategy:
\begin{itemize}
    \item Achieves service levels close to full-information policies.
    \item Incurs significantly lower communication costs.
    \item Outperforms heuristic and periodic sharing strategies.
\end{itemize}

The approach effectively optimizes the trade-off between communication and operational performance by sharing information only when its value exceeds its cost.

\section{Conclusion}\label{sec:conclusion}

We presented a novel application of decentralized control theory to a multi-retailer supply chain problem, formulating the joint design of communication and inventory control strategies. By reframing the problem in a common information framework, we derived a centralized POMDP whose solution yields optimal policies for when retailers should share data and how they should act on available information. The model captures the fundamental trade-off between improved operational decisions (through better information) and the incurred communication cost.

In a realistic setting with stochastic demand and inventory dynamics, we demonstrated that the optimal policy has an intuitive structure: it calls for communication only during significant demand changes or high uncertainty, and otherwise saves cost by withholding routine information. This adaptive approach outperforms both always-share and never-share extremes, as well as simple periodic sharing schedules, by maintaining high service levels with minimal communication.

Our work bridges ideas from supply chain management (which traditionally studies the value of information) and multi-agent decision-making (which provides tools for designing optimal strategies under asymmetric information). For practitioners, our results suggest that state-dependent communication policies (e.g., event-triggered data sharing among retail partners) can achieve most of the benefits of full information sharing at a much lower cost. This could inform the design of demand monitoring systems or collaborative agreements in retail consortia—for example, retailers could agree to share point-of-sale data only when sales deviate beyond certain thresholds, rather than constantly streaming all data.

There are several avenues for further research. First, scaling up the solution method to larger networks (more retailers, multi-echelon supply chains) is a key challenge—one could explore mean-field approximations or hierarchical decompositions to manage complexity. Second, extending the model to include lead times and backorders would increase realism; this introduces additional state variables (outstanding orders, backlogged demand) which the common information approach should handle (at the cost of larger state space). Third, considering strategic behavior and incentives is important for real-world adoption: in our cooperative model, all retailers had a common goal, but in practice independent firms might need incentive-compatible mechanisms (contracts) to follow the recommended sharing policy. Embedding our results into a game-theoretic framework (e.g., using mechanism design so that sharing decisions align with each retailer’s profit motive) would be an interesting extension. Lastly, one could integrate learning into the framework—for example, if demand model parameters are initially unknown, the communication strategy might also serve to coordinate learning (exploration) among retailers.

In conclusion, this paper demonstrates that optimal communication and control in multi-agent systems is a powerful paradigm for supply chain optimization. By judiciously sharing information, retailers can effectively operate as a coordinated team, reaping the benefits of collaboration while minimizing unnecessary communication. We believe the insights and methods developed here can be applied to other domains where decentralized actors must decide on both information exchange and operational actions (e.g., distributed energy grids, networked robotics, or multi-department inventory systems), leading to smarter and more efficient collaborative strategies.

% \section*{References}

\end{document}